\documentclass[12pt,letterpaper,reqno]{amsart}

\usepackage{times}
\usepackage[T1]{fontenc}
\usepackage{mathrsfs}
\usepackage{latexsym}
\usepackage[dvips]{graphics}
\usepackage{epsfig}
\usepackage{amsmath,amsfonts,amsthm,amssymb,amscd}
\input amssym.def
\input amssym.tex

\newcommand\be{\begin{equation}}
\newcommand\ee{\end{equation}}
\newcommand\bea{\begin{eqnarray}}
\newcommand\eea{\end{eqnarray}}
\newcommand\bi{\begin{itemize}}
\newcommand\ei{\end{itemize}}
\newcommand\ben{\begin{enumerate}}
\newcommand\een{\end{enumerate}}
\newcommand\bc{\begin{center}}
\newcommand\ec{\end{center}}
\newcommand\ba{\begin{array}}
\newcommand\ea{\end{array}}

\newtheorem{thm}{Theorem}[section]

\newtheorem{cor}[thm]{Corollary}
\newtheorem{lem}[thm]{Lemma}
\newtheorem{prop}[thm]{Proposition}

\theoremstyle{definition}
\newtheorem{rek}[thm]{Remark}

\begin{document}

\title{Finite groups with an automorphism cubing a large fraction of elements}

\author{Peter Hegarty}
\email{hegarty@math.chalmers.se} \address{Mathematical Sciences,
Chalmers University Of Technology and G\"oteborg University,
G\"oteborg, Sweden}


\date{\today}

\begin{abstract} We investigate the possible structures imposed on a 
finite group by its possession of an automorphism sending a large fraction of 
the group elements to their cubes, the philosophy being that this
should force the group to be, in some sense, close to abelian. We 
prove two main theorems. In the first, we completely classify all finite groups
with an automorphism cubing more than half their elements. All such 
groups are either nilpotent class 2 or possess an abelian subgroup of index 2.
For our second theorem, we show that if a group possesses an 
automorphism sending more than $4/15$ of its elements to their cubes, then it 
must be solvable. The group $A_{5}$ shows that this result is best-possible. 
\par Both our main findings closely parallel results of prevous authors on
finite groups possessing an automorphism which inverts many group elements.
The technicalities of the new proofs are somewhat more subtle, and also 
throw up a nice connection to a basic problem in combinatorial number 
theory, namely the study of subsets of finite cyclic groups which avoid
non-trivial solutions to one or more translation invariant linear equations.  
\end{abstract}


\maketitle

\setcounter{equation}{0}

\setcounter{equation}{0}

\section{Introduction}

Let $n$ be an integer. A group $G$ is said to be \emph{$n$-abelian} if 
the map $x \mapsto x^{n}$ is an endomorphism of $G$. It is a simple 
observation that, for $n = -1$ or $2$, an $n$-abelian group is 
abelian. The fact that there exist non-abelian groups of every 
exponent greater than or equal to three means that this observation does
not extend to any other value of $n$. However, Alperin \cite{Alp} obtained an 
elegant classification of $n$-abelian groups for every $n > 0$, his 
result being that a group is $n$-abelian if and only if it is a 
homomorphic image of a subgroup of the direct 
product of an abelian group, a group of exponent dividing $n$ and a 
group of exponent dividing $n-1$. In particular, for 
$n = 3$ this implies that a group for which the map $x \mapsto x^{3}$
is an injective endomorphism must also be abelian. 

Suppose $n \in \{-1,2,3\}$. For finite groups, the following questions
now arise \\ naturally :
\\
\\
1. Is there a constant $c_{n} < 1$ such that any finite group $G$ 
possessing an automorphism sending more than $c_{n} |G|$ elements to their
$n$:th powers is abelian ?
\\
\\
2. More generally, for what constants $c^{\prime}_{n} < 1$ can we produce 
an $\lq$elegant' (in some sense which is generally acceptable)
classification of finite groups $G$ possessing an automophism
sending more than $c^{\prime}_{n} |G|$ elements to their $n$:th powers ?
The groups appearing in the classification should all, in some sense, be 
$\lq$close' to abelian. 
\\
\\
Regarding Question 1, it is known that $c_{-1} = c_{3} = 3/4$ and $c_{2} = 
1/2$ : see \cite{Mil}, \cite{Mac1} and \cite{Zim} respectively. For
each prime $p$, let $\mathscr{G}_{p}$ denote the collection of finite
groups whose order is divisible by $p$ and by no smaller prime. 
Restricting attention to groups in $\mathscr{G}_{p}$ it is also 
known that $c_{n} = 1/p$ for each $n \in \{-1,2,3\}$ and for 
every odd $p$ : see \cite{LM2}, \cite{L} and \cite{Mac1}.  
\\
\\
Regarding Question 2, there is also a lot known. For each odd $p$, 
complete classifications are known of those
groups in $\mathscr{G}_{p}$ possessing an automorphism which
sends exactly $1/p$ of the group elements to their
inverses \cite{LM2}, squares \cite{L} respectively cubes \cite{DM}. For
even order groups there are the following results :
\\
\\
$n = -1$ : In what is probably the most significant paper in 
this area, Liebeck and MacHale \cite{LM1} provided a concise
classification of those groups admitting an automorphism which inverts more
than half their elements. MacHale and the author \cite{HM}
extended this classification to include groups admitting an 
automorphism which inverts exactly half the group elements, but already here
the classification is considerably more detailed. 
\\
\\
$n = -2$ : the author \cite{H}, improving upon results in
\cite{Zim}, classified neatly all even order groups 
possessing an automorphism squaring more than one-sixth of their elements. 
I also provided partial information at exactly one-sixth, but not 
a full classification. 
\\
\\
The missing piece in this jigsaw is a classification analogous to those 
above when $n = 3$. The main purpose of this paper is to provide this
missing piece (Theorem \ref{thm:mainclass} below). 
It is important to note here that all
the fractions appearing in these classifications (including ours) appear
to be optimal, i.e.: a reasonable corresponding description seems 
impossible for any smaller value of the fraction in question. In this
sense, we think that Theorem \ref{thm:mainclass} 
really does put a finishing touch to the
body of work outlined above. 
\\
\\
The methods introduced in \cite{LM1} provide the basis for much of the
subsequent investigations in the papers cited above. Let $n \in \{-1,2\}$. 
If an automorhpism $\alpha$ of a group $G$ sends a large fraction of the 
elements to their $n$:th powers, then for a large fraction of pairs $x,y$ of
elements the relation $x^{n}y^{n} = (xy)^{n}$ holds, and hence $[x,y] = 1$. 
Liebeck and MacHale exploit this information by focusing attention on a 
subgroup $H$ of $G$ of maximal order satisfying $h\alpha = h^{n} \; \forall
\; h \in H$, and considering the (right) coset decomposition of $H$ in $G$. 
If $x$ is any element of $G \backslash H$ such that $x\alpha = x^{n}$ then 
$\{h \in H : (hx)\alpha = (hx)^{n}\} = C_{H}(x)$ is a proper subgroup of 
$H$, by definition of the latter. These observations form the basis of a
counting argument which eventually leads to the kinds of 
results we refer to above. 
\\
\\
For $n = 3$ we want to apply the same type of argument, but we run into an 
immediate difficulty, namely : the relation $x^{3}y^{3} = (xy)^{3}$ 
does not on its own imply that $x$ and $y$ commute. The main 
contribution of the present paper is to remove this obstacle to 
obtaining results for $n = 3$ which are as good as those for $n \in \{-1,2\}$. 
The technical results obtained in Section 2 for this purpose 
are thus, in my opinion, the 
real heart of the paper, especially since they establish an unexpected 
connection between our problem and a fundamental problem in 
combinatorial number theory, namely the study of sets of integers which 
contain no non-trivial solutions to one or more translation invariant 
linear equations. These connections, which may be of independent
interest, are summarised in Proposition \ref{prop:numbertheory} below.
\\
\\
The final classification obtained in Theorem \ref{thm:mainclass} is 
almost identical to the one in \cite{LM1}, except for 
obvious extra conditions on the 3-part of $G$. This is, in fact, not a 
surprise, once the machinery in Section 2 has been developed, though the path
to the final result is still more difficult than in \cite{LM1}. Section 3 is 
devoted to the proof of this theorem. To illustrate further the 
effectiveness of our machinery, we devote Section 4 to a proof of the fact 
(Theorem \ref{thm:solvable}) that a finite group admitting an 
automorphism sending more than $4/15$:ths of its elements to their
cubes must be solvable. This mirrors analogous results for inverses 
\cite{P} and squares \cite{H}, where the corresponding constants are $4/15$ and
$7/60$ respectively. Curiously the same group, namely $A_{5}$, illustrates that
all three constants are optimal.         
\\
\\
The final section (Section 5) provides a brief summary of our findings and 
a discussion of outstanding issues.    

\setcounter{equation}{0}
\section{Preliminary Lemmas and Connections to Number Theory}\label{sec:preliminaries}

First let us fix some notation. If $G$ is a finite group and $\alpha$ an 
automorphism of $G$, we denote
\begin{eqnarray*}
T_{3,\alpha} := \{g \in G : g\alpha = g^{3} \}
\end{eqnarray*}
and
\begin{eqnarray*}
r_{3}(G,\alpha) := {|T_{3,\alpha}| \over |G|}.
\end{eqnarray*}
If $N$ is an $\alpha$-invariant subgroup of $G$, we denote by $\alpha_{N}$
the restriction of $\alpha$ to $N$. If, in addition, $N \lhd G$ then
the induced automorphism of $G/N$ is denoted $\alpha^{N}$. 
\par We reserve the letter $H$ for a subgroup of $G$ contained inside
$T_{3,\alpha}$. In Section 3, but not otherwise, 
we will further reserve $H$ to denote a subgroup of maximum order 
with this property. For $x,y \in G$, the commutator 
$x^{-1}y^{-1}xy$ is denoted $[x,y]$. Finally, for $n > 0$, the cyclic 
group of order $n$ is denoted $\mathbb{Z}_n$.  
\\
\\
In the following lemmas, consider a group $G$ and an automorphism $\alpha$ as 
given. The proofs of the first two results are obvious :

\begin{lem}\label{lem:induced} If $N \lhd G$ is $\alpha$-invariant, then 
$r_{3}(G,\alpha) \leq r_{3}(G/N,\alpha^{N})$. 
\end{lem}

\begin{lem}\label{lem:three} If $x \in T_{3,\alpha}$ then $C_{G}(x) = 
C_{G}(x^{3})$. In particular, if $H \subseteq T_{3,\alpha}$ then 
$(H:C_{H}(x))$ is not divisible by three.
\end{lem}

The next two results are also easy :

\begin{lem}\label{lem:eltwoab} If $H \subseteq T_{3,\alpha}$, $x \in 
T_{3,\alpha}$ and $H/C_{H}(x^{2})$ is elementary 2-abelian, then 
$hx \in T_{3,\alpha} \Leftrightarrow [h,x] = 1$. 
\end{lem}

\begin{proof}
Suppose $hx \in T_{3,\alpha}$. Then 
\begin{eqnarray*}
(hx)\alpha = (hx)^{3} = h\alpha x \alpha = h^{3} x^{3},
\end{eqnarray*}
which implies that $h^{2}x^{2} = (xh)^{2}$. But our 
assumptions imply that $[h^{2},x^{2}] = 1$, thus $(xh)^{2} = 
h^{2}x^{2} = x^{2} h^{2}$, from which it follows that $[h,x] = 1$. 
\end{proof}

\begin{lem}\label{lem:abba} Suppose each of $a,b,ab$ and $ba$ is in 
$T_{3,\alpha}$. Then $[a,b] = 1$.
\end{lem}

\begin{proof}
As in the proof of the previous lemma, we can deduce immediately from 
our assumptions that 
\begin{eqnarray*}
a^2b^2 = (ba)^2, \;\;\;\; b^2a^2 = (ab)^2.
\end{eqnarray*}
But then $a^2b^3 = (ba)^2b = b(ab)^2 = b^3 a^2$, so $[a^2,b^3] = 1$. But 
then $[a^2,b] = 1$ by Lemma \ref{lem:three}, so now 
$(ba)^2 = a^2b^2 = b^2a^2$, thus $[a,b] = 1$. 
\end{proof}

The next result is the crucial one :

\begin{lem}\label{lem:ap} Suppose each of $a,b,ab$ and $a^{-1}b$ is in 
$T_{3,\alpha}$. Then $[a,b] = 1$.
\end{lem}

\begin{proof}
As previously, we can deduce immediately from our
assumptions that 
\be
a^{2}b^{2} = (ba)^{2}
\ee
and 
\be
a^{-2}b^2 = (ba^{-1})^2.
\ee
From these and the identity 
\begin{eqnarray*}
[x,yxy^{-1}] = [x^{-1}y)^3 (y^{-1}xy^{-1})^3 (ya^{-1}yxy^{-1})^3
\end{eqnarray*}
it is easily deduced that $bab^{-1} \in C_{G}(a)$, from which we also
deduce, using (2.1), that $a^{2}ba \in C_{G}(b)$. Thus 
\be
a^{2}bab^{-1} \in C_{G}(a) \cap C_{G}(b) \supseteq C_{G}(ab).
\ee
Now, since $G$ is finite, there exists a positive integer $n$ such that
$a^{n} \in C_{G}(b)$. First suppose $n$ is even, say $n = 2k$. Then, by (2.3),  
$(a^{2}bab^{-1})^k = a^{2k}bab^{-1} \in C_{G}(b)$, hence $a^{k} \in 
C_{G}(b)$. Thus we may in fact assume $n$ is odd, say $n = 2k+1$. 
\par Then, using (2.3) again, we have that $(a^{2}bab^{-1})^{k+1} = 
a^{2k+1}(ab)a^{2k+1}a^{-k}b^{-1} \in C_{G}(ab)$, which implies that 
$a^{-k}b^{-1} \in C_{G}(ab)$ and hence that $(ba^{k})^3 \in C_{G}(ab)$.
\par But (2.3) also implies that $b^{-1}a^2b = abab^{-1}$, hence that
$b^{-1}a^{2k}b = a^kba^kb^{-1}$, and in turn that
$a^{2k}b^2 = (ba^k)^2$. Thus 
\be
(ba^k)^3 = a^{2k} b^3 a^k \in C_{G}(ab).
\ee
Furthermore, by Lemma \ref{lem:three} we may assume that $n$ is not divisible
by three, so that $k = 3l$ or $k = 3l+2$ for some $l$. 
\par First suppose $k = 3l$. Then (2.4) says that $a^{6l}b^3 a^{3l} \in
C_{G}(ab)$. But $a^{6l}b^3a^{3l} = (a^{2l}ba^l)\alpha$ and $ab \in T_{3,\alpha}$, 
hence 
\be
a^{2l}ba^l \in C_{G}(ab), 
\ee
by Lemma \ref{lem:three}. But, going back to (2.3), we have $(a^2bab^{-1})^l =
a^{2l}ba^lb^{-1} \in C_{G}(ab)$ which, together with (2.5), implies that
$b^{-1} \in C_{G}(ab)$, hence that $[a,b] = 1$ as required. 
\par Alternatively, if $k = 3l+2$, then $n \equiv -1 \; ({\hbox{mod $6$}})$
so $n^2 \equiv 1 \; ({\hbox{mod $6$}})$. Thus if we work with $n^2$ instead
of $n$ we will get the same conclusion, namely that $[a,b] = 1$, and so the 
lemma is proved.
\end{proof}     

\begin{rek}
In the above proof we have used the finiteness of $G$ to guarantee that 
some power of $a$ commutes with $b$. Hence the proof goes through in any
torsion group, for example. But we do not know whether these
restrictions are really necessary, or whether the lemma holds in arbitrary
groups.
\end{rek}

\begin{cor}\label{cor:ap2} Suppose each of $a,b,ab$ and $a^{-2}b$ is in 
$T_{3,\alpha}$. Then $[a,b]=1$.
\end{cor}

\begin{proof}
The assumptions imply that 
\begin{eqnarray*}
a^2b^2 = (ba)^2, \;\;\;\; a^{-4}b^2 = (ba^{-2})^2.
\end{eqnarray*}
Then 
\begin{eqnarray*}
a^6 = (a^2b^2)(a^{-4}b^2)^{-1} = baba^3b^{-1}a^2b^{-1}, 
\end{eqnarray*}
and so 
\begin{eqnarray*}
a^6b^3 = baba^3b^{-1}(a^2b^2) = (ba)^2a^3ba = a^2b^2a^3ba, 
\end{eqnarray*}
from which it follows that $b^{-2}a \in C_{G}(a^3b^3)$. Then Lemma 
\ref{lem:three} implies that, in fact, $b^{-2}a \in C_{G}(ab)$. But then
$(b^{-2}a)(ab) = (ab)(b^{-2}a)$, hence $a^{-2}b^2 = (ba^{-1})^2$, which
implies that $a^{-1}b \in T_{3,\alpha}$. Now the result follows
from Lemma \ref{lem:ap}.
\end{proof}

\begin{rek}
Another corollary of Lemma \ref{lem:ap} is that if $a,b,ab$ and $a^2b$ are all
in $T_{3,\alpha}$, then $[a,b] = 1$. This follows immediately from the 
lemma upon making the variable substitutions $a^{\prime} := a$, $b^{\prime} := 
ab$. Similarly, if $\{a,b,ab,a^3 b\} \subseteq T_{3,\alpha}$ then $[a,b]=1$.
This follows from Corollary \ref{cor:ap2} upon substituting
$a^{\prime} := a^{-1}$, $b^{\prime} := ab$.
We do not know if it is possible to obtain further results like these. 
One may ask : does there exist any integer $n \not\in \{-1,\pm2,3\}$ such that,
if $\{a,b,ab,a^n b\} \subseteq T_{3,\alpha}$ then one must have 
$[a,b] = 1$ ? We suspect that there are no other such $n$.  
\end{rek}

Let $H$ be a subgroup of $G$ such that $H \subseteq T_{3,\alpha}$. Thus $H$
is abelian. Let $x \in T_{3,\alpha}$. Then clearly, $\{h \in H : hx \in
T_{3,\alpha}\}$ consists of entire cosets in $H$ of $C_{H}(x)$. Thus
the set $Hx \cap T_{3,\alpha}$ may be identified with a subset, which we denote
$\mathscr{T}(H,x)$, of the abelian group $H/C_{H}(x)$. The last two
results now immediately yield the following, which establishes the 
connection referred to earlier 
between our work and combinatorial number theory :

\begin{prop}\label{prop:numbertheory} For any subgroup $H \subseteq 
T_{3,\alpha}$ and any $x \in T_{3,\alpha}$, the subset $\mathscr{T}(H,x)$ of 
the abelian group $H/C_{H}(x)$, written additively, contains no 
non-trivial solutions to either of the translation invariant linear equations 
$a+b=2c$, $a+2b=3c$. In particular, it contains no 3-term arithmetic
progressions. 
\end{prop}

\begin{proof}
This follows directly from Lemma \ref{lem:ap} and Corollary \ref{cor:ap2}. Note
that a 3-term arithmetic progression is just a solution
to $a+b=2c$ with $a \neq b$ (we allow $a=c$, which can arise in groups of even 
order). 
\end{proof}

Let $f(x_{1},...,x_{n}) \in \mathbb{Z}[x_{1},...,x_{n}]$ 
be any translation invariant linear 
function, i.e.: $f(x_{1},...,x_{n}) = \sum_{i=1}^{n} a_{i} x_{i}$ where 
$a_{i} \in \mathbb{Z}$ and $\sum_{i=1}^{n} a_{i} = 0$. It is known that if 
$A \subseteq \mathbb{N}$ has non-zero upper asymptotic density then 
$A$ must contain a non-trivial solution to $f(x_{1},...,x_{n}) = 0$. This is 
an easy consequence of the celebrated theorem of Szemer\'{e}di stating
that if $A \subseteq \mathbb{N}$ has non-zero upper asymptotic density, then 
$A$ contains arbitrarily long arithmetic progressions. For a 
discussion of these results, inlcuding a formal definition of what is meant
by a $\lq$non-trivial solution' of a translation invariant linear
equation, see \cite{R}. Note that, for an equation in three variables, like
those appearing in Proposition \ref{prop:numbertheory},
non-trivial means simply that $x_{1},x_{2},x_{3}$ are not all equal. 
\par These results have immediate corollaries in finite cyclic groups, namely, 
as $n \rightarrow \infty$, if $A \subseteq \mathbb{Z}_n$ 
contains no non-trivial solutions to $f(x_{1},...,x_{n}) \equiv 0 
\; ({\hbox{mod $n$}})$, then $|A| = o(n)$. This is, in fact, what we will 
use in Section 4 of this paper, where the subgroup $H$ will always be a
cyclic group generated by a single element of $T_{3,\alpha}$. It is worth noting
though, that corresponding results exist for arbitrary finite abelian 
groups : for an up-to-date treatment of these matters, see for example 
\cite{GT}.
\\
\\
Speaking somewhat 
informally, Proposition \ref{prop:numbertheory} and the above
results from number theory imply the following : Let $G$ be a finite
group possessing an automorphism $\alpha$ for which $r_{3}(G,\alpha)$ is large.
Then either there is a correspondingly large proportion of 
commuting pairs of group elements after all (as would be the case if
we replaced $r_{3}$ by $r_{-1}$ or $r_{2}$), or most of the elements of 
$T_{3,\alpha}$ have small order.


\setcounter{equation}{0}
\section{Proof of Classification Theorem}\label{sec:theoremA}

The purpose of this section is to prove the following theorem :

\begin{thm}\label{thm:mainclass} The finite group $G$ admits an 
automorphism $\alpha$ for which $r_{3}(G,\alpha) > 1/2$ if and only if $G$
has one of the following structures :
\\
\\
\textbf{I.} $G$ is abelian and $(|G|,3) = 1$. 
\\
\\
\textbf{II.} $G$ is non-abelian with a 
normal Sylow 3-subgroup $S$ satisfying the following
conditions :
\par (a) $S \subseteq K$ where $(G:K) = 2$ and $K$ is abelian,
\par (b) $S \cap Z(G) = \{1\}$. 
\\
In particular, if $(|G|,3) = 1$ then it suffices for $G$ to have an abelian
subgroup of index 2. 
\\
\\
\textbf{III.} $G$ is nilpotent class two and $(|G|,3) = 1$. All 
Sylow $p$-subgroups, for $p > 2$, are abelian. The Sylow 2-subgroup $S_{2}$ 
has one of the following structures :
\\
\\
\textbf{(i)} $S_{2}^{\prime} \cong C_{2} = \; <z>$, say. 
$S_{2}/Z(S_{2})$ is elementary 
abelian, generated by $Zx_{1},...,Zx_{k}$, $Za_{1},...,Za_{k}$, subject to the
following commutator relations :
\begin{eqnarray*}
[x_{i},x_{j}]=[a_{i},a_{j}]=[a_{i},x_{j}]=1 \;\; {\hbox{whenever $i \neq j$}}, 
\;\; [a_{i},x_{i}] = z.
\end{eqnarray*}
\textbf{(ii)} $S_{2}^{\prime} \cong C_{2} \times C_{2} = \; 
<z_{1}> \times <z_{2}>$,
say. $S_{2}/Z(S_{2})$ is elementary abelian of order 16, generated by
$Zx_{1},Zx_{2},Za_{1},Za_{2}$, subject to the following commutator relations :
\begin{eqnarray*}
[x_{i},x_{j}]=[a_{i},a_{j}]=[a_{i},x_{j}]=1 \;\; {\hbox{whenever $i \neq j$}}, 
\;\; [a_{i},x_{i}] = z_{i}.
\end{eqnarray*}
\end{thm}
First let us deal with the $\lq$if' part of the theorem by constructing an
explicit automorphism $\alpha$ of each type of group such that
$r_{3}(G,\alpha) > 1/2$.
\\
\\
\textbf{I.} The map $\alpha : g \mapsto g^3 \; \forall \; g \in G$ is an
automorphism and $r_{3}(G,\alpha) = 1$. 
\\
\\
\textbf{II.} If $x \in G\backslash K$ then $(|x|,3) = 1$ since, if 
$x^{3^m n} = 1$ then, by normality of $S$ and commutativity of $K$, 
we have $x^n \in S \cap Z(G) = \{1\}$. Now fix any choice of $x \in G\backslash
K$ and define the map $\alpha : G \rightarrow G$ as follows :
\begin{eqnarray*}
k \alpha := k^2 x^{-1}kx \; \forall \; k \in K, \;\;\; x\alpha := x^3, \;\;\;
(kx)\alpha := k\alpha x\alpha \; \forall \; k \in K.
\end{eqnarray*}
It is easily checked that $\alpha$ is well-defined and thus a homomorphism.
Furthermore, $\alpha$ is one-to-one on $K$ since
$k^2 x^{-1}kx = 1 \Leftrightarrow x^3 = (xk^{-1})^3 \Leftrightarrow
x = xk^{-1} \Leftrightarrow k = 1$, where we have used the fact that
$(|g|,3)=1$ for all $g \in G\backslash K$. Thus $\alpha \in$ Aut$(G)$.
Finally, it is also easily verified that $T_{3,\alpha} = Kx \sqcup
C_{K}(x)$, hence $r_{3}(G,\alpha) = {n+1 \over 2n} > {1 \over 2}$, where
$(K:C_{K}(x)) = n$. 
\\
\\
\textbf{III.} Let $A$ be the abelian subgroup of $G$ generated by 
$Z(G)$ and $a_{1},...,a_{k}$. The map $\alpha : G \rightarrow G$ defined by 
\begin{eqnarray*}
(ax_{1}^{\epsilon_{1}}x_{2}^{\epsilon_{2}} \cdots x_{k}^{\epsilon_{k}})\alpha := 
a^3 x_{1}^{3\epsilon_{1}}x_{2}^{3\epsilon_{2}} \cdots x_{k}^{3\epsilon_{k}} 
\;\;\; \forall \; a \in A, \; \epsilon_{i} \in \{0,1\}, i = 1,...,k,
\end{eqnarray*}
is easily seen to be an automorphism of $G$ such that $r_{3}(G,\alpha) = 
{2^k + 1 \over 2^{k+1}}$. In particular, for groups of 
type \textbf{(ii)} we have $r_{3}(G,\alpha) = 5/8$. 
For more details, see \cite{LM1}.

\begin{rek}
For each of the groups $G$ in Theorem \ref{thm:mainclass}, 
it is easy to show that there is no 
$\beta \in$ Aut$(G)$ such that $r_{3}(G,\beta) > r_{3}(G,\alpha)$, where 
$\alpha$ is the automorphism constructed above. See \cite{LM1} for similar
remarks. 
\end{rek}

Now we turn to the $\lq$only if' part of the theorem. Fix a group $G$ and an 
automorphism $\alpha$ for which $r_{3}(G,\alpha) > 1/2$. For the remainder
of this section, $H$ will denote a subgroup of $G$ of maximum order 
subject to the condition that $H \subseteq T_{3,\alpha}$. The center of 
$G$ will be denoted simply by $Z$.  

\begin{lem}\label{lem:center} $H \supseteq Z$.
\end{lem}

\begin{proof}
By considering a decomposition of $G$ into cosets of $Z$ we see 
that if $r_{3}(G,\alpha) > 1/2$ then $Z \subseteq T_{3,\alpha}$. Since 
$<Z,x> \; \subseteq T_{3,\alpha}$ for any $x \in T_{3,\alpha}$, maximality of 
$H$ implies that $H \supseteq Z$. 
\end{proof}

In the notation of Proposition \ref{prop:numbertheory} let us denote
\begin{eqnarray*}
t(H,x) := {|\mathscr{T}(H,x)| \over |H/C_{H}(x)|}.
\end{eqnarray*}
In this section we only need some very weak consequences of the machinery
developed in Section 2, namely :

\begin{lem}\label{lem:cosets} Let $x \in T_{3,\alpha} \backslash H$. Then 
$|Hx \cap T_{3,\alpha}| \leq \frac{1}{2} |H|$. Hence every 
right-coset of $H$ in $G$ intersects $T_{3,\alpha}$. In particular, $H$
is not properly contained in any other abelian subgroup of $G$. 
Moreover, if $Hx \neq 
Hx^{-1}$ and $(H:C_{H}(x)) = n > 2$ then $|Hx \cap T_{3,\alpha}| + 
|Hx^{-1} \cap T_{3,\alpha}| \leq \frac{2}{n} \left( 1 + \frac{1}{4}
\lfloor n-1 \rfloor \right) |H|$. 
\end{lem}

\begin{proof}
By maximality of $H$, the group $H/C_{H}(x)$ must be non-trivial. Then it is 
an elementary consequence of Proposition \ref{prop:numbertheory} that 
$t(H,x) \leq 1/2$. This implies the first assertion of the lemma. The second 
one follows immediately and then the third from the definition of $H$.  
\par For the final assertion, let $K := C_{H}(x) = C_{H}(x^{-1})$ and 
consider $H/K$ as an additive group. Let 
\begin{eqnarray*}
S_{+} := \mathscr{T}(H,x) \backslash \{0\}, \;\;\;\; 
S_{-} := \mathscr{T}(H,x^{-1}) \backslash \{0\}.
\end{eqnarray*}
Lemma \ref{lem:ap} implies that 
\begin{eqnarray*}
S_{+} \cap (-S_{+}) = S_{-} \cap (-S_{-}) = S_{+} \cap S_{-} = \phi,
\end{eqnarray*}
from which the result follows.
\end{proof}
Let $(G:H) = m$ and 
\be
G = H \sqcup Hx_{2} \sqcup \cdots \sqcup Hx_{m}
\ee
be a right-coset decomposition of $H$ in $G$ such that 
$x_{i} \in T_{3,\alpha}$ for every $i \in \{2,...,m\}$. Such a 
decomposition exists by 
Lemma \ref{lem:cosets}. Then 
\be
r_{3}(G,\alpha) = {1 \over m} \left( 1 + \sum_{i=2}^{m} t(H,x_{i})
\right).
\ee
The next result will set us up nicely for the remainder of the proof of our
theorem :

\begin{lem}\label{lem:strange} Assuming $r_{3}(G,\alpha) > 1/2$ we must have
$r_{3}(G,\alpha) = {n+1 \over 2n}$ for some $n \in \mathbb{N}$. Moreover, 
in a right-coset decomposition of $H$ in $G$ as in (3.1), we must have
\\ $(H:C_{H}(x_{i})) > 2$ for at most one index $i$ and $\mathscr{T}(H,x_{i}) = 
\{0\}$ for every index $i$. 
\end{lem}

\begin{proof}
Let there be $k$ indices $i$ for which $(H:C_{H}(x_{i})) > 2$.
\\
\\
{\sc Case I} : $k = 0$. 
\\
\\
In this case, by (3.2), $r_{3}(G,\alpha) = {m+1 \over 2m}$ where
$(G:H) = m$. Clearly, $\mathscr{T}(H,x_{i}) = \{0\}$ whenever 
$(H:C_{H}(x_{i})) = 2$. 
\\
\\
{\sc Case II} : $k = 1$.
\\
\\
Suppose $(H:C_{H}(x_{m})) = n > 2$. If $(G:H) = 2$ then Lemma 
\ref{lem:eltwoab} implies that $\mathscr{T}(H,x_{2}) = \{0\}$ and so
$r_{3}(G,\alpha) = {n+1 \over 2n}$. Otherwise we must have 
$Hx_{m} = Hx_{i}x_{j}$ for some $i,j < m$ and so $H/C_{H}(x_{m}) 
\cong C_{2} \times C_{2}$. Thus $n = 4$, $\mathscr{T}(H,x_{m}) = \{0\}$ and
$r_{3}(G,\alpha) = {2m+1 \over 4m}$. 
\\
\\
{\sc Case III} : $k = 2$.
\\
\\
Let $i,j$ be the indices such that $(H:C_{H}(x_{i})) = n_{i} > 2$ and 
$(H:C_{H}(x_{j})) = n_{j} > 2$. By Lemma \ref{lem:three}, in fact  $n_{i}
\geq 4$ and $n_{j} \geq 4$. If neither $Hx_{i}^2 = Hx_{j}$ nor $Hx_{j}^2 =
Hx_{i}$ holds, then Lemma \ref{lem:eltwoab} implies that
$\mathscr{T}(H,x_i) = \mathscr{T}(H,x_j) = \{0\}$ and (3.2) gives
$r_{3}(G,\alpha) \leq 1/2$, a contradiction. Hence, we may assume that 
$Hx_i^2 = Hx_j$, say. But then, by
Lemma \ref{lem:three}, there is a third coset, namely $Hx_i^3$, such that
$(H:C_{H}(x_i^3)) > 2$. Thus {\sc Case III} cannot arise.
\\
\\
{\sc Case IV} : $k > 2$.
\\
\\
Let $Hy=Hy_1,Hy_2,...,Hy_k$ be a complete set of cosets of $H$ in $G$ for
which $y_i \in T_{3,\alpha}$ and $(H:C_{H}(y_i)) = n_i \geq 4$. If $y_i^2
\in H$ then $t(H,y_i) = 1/n_i$ by Lemma \ref{lem:eltwoab}, so if this were
the case for every $i = 1,...,k$ then (3.2) would imply that
$r_3(G,\alpha) < 1/2$.
\par Without loss of generality, suppose $y^2 \not\in H$. Thus the cosets
$Hy$ and $Hy^{-1}$ are distinct. If $(H:C_{H}(y^2))=2$ then Lemma
\ref{lem:eltwoab} and (3.2) again give the contradiction that
$r_3(G,\alpha) \leq 1/2$. Thus $Hy^2 = Hy_j$ for some $j$. But, using both
Lemmas \ref{lem:eltwoab} and \ref{lem:cosets} this time, we'll get
the same contradiction if $(H:C_{H}(y^4)) \leq 2$. In particular, we may
assume that $y^4 \not\in H$ and hence that the four cosets $Hy, Hy^{-1},
Hy^2, Hy^{-2}$ are distinct. Grouping these in two pairs and using 
Lemma \ref{lem:cosets} again, we arrive at the same contradiction unless
$(H:C_{H}(y))=5$ and $t(H,y) = 2/5$. 
In this case, maximality of $H$ means that 
$y^5 \in H$. But then we claim that, in fact, $\mathscr{T}(H,x) = \{0\}$. For
if $hy \in T_{3,\alpha}$ then so is $(hy)\alpha^{2} = h^9y^9$, and hence
$h^{-1}y^{-1} \in T_{3,\alpha}$. But then $h \in C_{H}(y)$ by Lemma
\ref{lem:abba}.
\par Thus {\sc Case IV} cannot arise either, and so the proof of 
Lemma \ref{lem:strange} is complete.  
\end{proof}

Let us call a right coset $Hx$ \emph{exceptional} if $(H:C_{H}(x)) > 2$. 
By Lemma \ref{lem:strange} there is at most one exceptional coset of 
$H$ in $G$. Moreover, we have 

\begin{cor}\label{cor:strange} Suppose $(G:H) > 2$. Then $h^2 \in Z$
for all $h \in H$. In fact, $x^2 \in Z$ whenever $x \in T_{3,\alpha}$ and
the coset $Hx$ is not exceptional. If $x \in T_{3,\alpha}$ and $Hx$
is exceptional, then $H/C_{H}(x) \cong C_2 \times C_2$, 
$x^2 \in H$ and $x^4 \in Z$. 
\end{cor}

\begin{proof}
Lemma \ref{lem:strange} immediately implies that $h^2 \in Z$ for all
$h \in H$. If $x \in T_{3,\alpha}$ and the coset $Hx$ is not exceptional, then
the subgroup $<C_{H}(x),x>$ has the same properties as $H$, so applying the 
lemma to it instead yields that $x^2 \in Z$. Suppose $Hx$ is exceptional. If
$x^2 \in H$ then $x^4 \in Z$, so suppose $x^2 \not\in H$. Then 
the subgroup $<C_{H}(x),x>$ has the same properties as $H$, and so 
$x^2 \in Z$, a contradiction. 
\end{proof}

Note that if $(G:H) = 2$ then $G$ is of type \textbf{I} or
\textbf{II} in Theorem \ref{thm:mainclass}. So henceforth we
shall always assume that $(G:H) > 2$. 
We require two further preparatory results before presenting the main 
body of our argument.

\begin{lem}\label{lem:normal} Suppose that for every possible choice of 
the subgroup $H$ we have that $H \lhd G$. Then there is an automorphism 
$\beta$ of $G$, possibly different from $\alpha$, such that 
$r_{3}(G,\beta) > 1/2$ and for which any corresponding $H_{\beta}$ is an 
abelian subgroup of maximum order in $G$. Moreover, $G$ is of type 
\textbf{III} in Theorem \ref{thm:mainclass}.
\end{lem}

\begin{proof}
From Corollary \ref{cor:strange} we know that $x^2 \in H$ for all $x \in 
T_{3,\alpha}$. If $H \lhd G$ this implies that $g^2 \in H$ for all $g \in G$. 
If the same is true for any possible choice of $H$ then, by Lemma
\ref{lem:strange}, it follows that $g^2 \in Z$ for all $g \in G$, since 
$Z$ is just the intersection of all the possible choices for $H$. 
\par Now let $x \in T_{3,\alpha}$ and $I_{x}$ be the inner automorphism of
$G$ which sends $g$ to $x^{-1}gx$. Since $g^2 \in Z$ for all $g$, it is easily
checked that $gx \in T_{3,\alpha}$ if and only if $g \in T_{3,I_{x}\alpha}$. 
Thus $r_{3}(G,\alpha) = r_{3}(G,I_x \alpha)$ for any $x \in T_{3,\alpha}$.
\par Now let $A$ be an abelian subgroup of maximum order in $G$. Since 
$r_3(G,\alpha) > 1/2$, there is some coset $Ax$  such that $x \in T_{3,\alpha}$
and $|Ax \cap T_{3,\alpha}| > \frac{1}{2} |A|$. But then 
$|A \cap T_{3,I_x \alpha}| > \frac{1}{2} |A|$, so $A \subseteq T_{3,\alpha}$
since $A$ is abelian. So we choose $\beta := I_x \alpha$. It remains to show
that $G$ is of type \textbf{III} in Theorem \ref{thm:mainclass}. This is 
highly non-trivial, but the argument parallels entirely that in 
Section 4 of \cite{LM1}, with very minor modifications. We thus
omit further details.   
\end{proof}

\begin{lem}\label{lem:indextwo} Suppose that $(H:Z) = 2$. Let $K := Z \cup 
G\backslash T_{3,\alpha}$. Then $K$ is an abelian subgroup of index 2 in 
$G$. 
\end{lem}

\begin{proof} The assumption implies that there is no exceptional coset, and
hence $x^2 \in Z$ for all $x \in T_{3,\alpha}$, by Corollary \ref{cor:strange}.
Thus if $a,b$ and $b^{-1}a$ are each in $T_{3,\alpha}$ then so is 
$b^2(b^{-1}a) = ba$, and so $[a,b] = 1$ by Lemma \ref{lem:ap}. By maximality of 
$H$, it follows that, for any $x \in T_{3,\alpha} \backslash Z$, we have 
$C_{G}(x) \cap T_{3,\alpha} = \; <Z,x>$. 
\par To show that $K$ is closed under multiplication, it suffices to show 
that if $g_1, g_2 \in K$ then $g_2^{-1}g_1 \in K$. Clearly this is the 
case if either $g_1$ or $g_2$ lies in $Z$. So suppose $\{g_1,g_2\}
\subseteq K \backslash Z$. Let $H = \; <Z,h>$. 
By Lemma \ref{lem:strange}, there
exist $x_1, x_2 \in T_{3,\alpha}\backslash Z$ such that $g_i = hx_i$ for 
$i = 1,2$. Then $g_2^{-1}g_1 = x_{2}^{-1}x_1$, and by the above observations,
this lies in $T_{3,\alpha}$ if and only if $[x_1,x_2] = 1$, hence if and only
if $x_2 \in \; <Z,x_1>$. But this will imply that either $g_2^{-1}g_1 \in Z$,
which is okay, or that $g_2 \in H\backslash Z$, 
contradicting that $g_2 \in K$. 
\par This proves that $K$ is closed, hence a subgroup of $G$. Clearly $(G:K) = 
2$ and, by its definition, we can write $G = K \sqcup Kx$, where 
$Kx \subset T_{3,\alpha}$. Then for any $k \in K$ we have that 
\begin{eqnarray*}
(kx)\alpha = (kx)^3 = k\alpha x\alpha = k\alpha x^3, 
\end{eqnarray*}
hence $k\alpha = kx^{-1}kxk$, since $x^2 \in Z$. But since this holds for any
choice of $x$ and $k$, it follows that $K$ is abelian.
\end{proof}

By Lemma \ref{lem:induced} and Corollary \ref{cor:strange} the 
induced automorphism $\alpha^{Z}$ of $G/Z$ sends more than half its 
elements to their inverses. By the main result of 
\cite{LM1} there are the following three possibilities :
\\
\\
\textbf{(A)} $G/Z$ is abelian.
\\
\\
\textbf{(B)} $G/Z$ is nilpotent class two with $(G/Z)^{\prime} \cong
C_2$ or $C_2 \times C_2$, and various other conditions.
\\
\\    
\textbf{(C)} $G/Z$ has an abelian subgroup of index 2. 
\\
\\
If \textbf{(A)} holds then we are done, by Lemmas \ref{lem:center} and 
\ref{lem:normal}. Next we deal with \textbf{(B)} by proving 

\begin{lem}\label{lem:factorb} Let $G$ be a group possessing an 
automorphism $\alpha$ for which $r_3(G,\alpha) > 1/2$. 
Suppose that $G$ is nilpotent of 
class at most 3 and that $(G/Z)^{\prime}$ is elementary abelian of order at
most 4. Then unless $G$ has an abelian subgroup of index 2, the class of $G$
is at most 2.
\end{lem}

Note that this will indeed deal with \textbf{(B)}, by Lemma 
\ref{lem:normal}.       

\begin{proof}
We consider a minimal counterexample to the lemma and obtain a contradiction. 
By the results in \cite{DM} we know that all Sylow $p$-subgroups of $G$, for 
$p > 2$, are abelian, so we may assume $G$ to be a $2$-group. Further, by
Lemma \ref{lem:normal}, we may assume that there is a choice of the 
subgroup $H$ which is not normal in $G$. We fix such a choice once and for
all. In the body of the text to follow, we shall assume that there are
no exceptional right cosets of $H$ in $G$. Some additional
technicalities arise otherwise, and these will be indicated by means of 
footnotes. 

Let $N := N_G(H)$. 
Since $G$ is nilpotent, we have a strict containment $H \subset
N$. We consider three cases :
\\
\\
{\sc Case 1} : $N$ contains an abelian subgroup of index 2, but 
$(N:H) > 2$.
\\
{\sc Case 2} : $(N:H) = 2$.
\\
{\sc Case 3} : $N$ contains no abelian subgroup of index 2. 
\\
\\
First consider {\sc Case 1}. Let $K$ denote the abelian subgroup of index 2. 
By Lemma \ref{lem:cosets}, $K$ does not contain $H$, so 
$(H:K\cap H) = 2$. But $K\cap H = Z(N)$. 
Since $N$ is $\alpha$-invariant, we can now apply Lemma \ref{lem:indextwo}
to it to conclude that it possesses an abelian subgroup $L$ of index 2, 
possibly different from $K$. Indeed, $L = (K\cap H) \cup 
N\backslash T_{3,\alpha}$. 
\par Suppose $L \lhd G$. Let $x \in T_{3,\alpha} \backslash N$ and 
$h \in Z(N)$. Then $x^{-1}hx \in L$. But $x^{-1}hx \in 
T_{3,\alpha}$ since $x^2 \in H$ (Corollary \ref{cor:strange}), thus 
$x^{-1}hx \in L \cap T_{3,\alpha} \subset H$. But 
$x \not\in N$ so if $(H:C_{H}(x)) = 2$, then $h \in 
C_{H}(x)$. Since
$x$ was chosen arbitrarily and there is at most one exceptional coset, 
it follows that $h \in Z$. Thus 
$(H:Z) = 2$ and so $G$ possesses an abelian subgroup of index 2.
\par So we may assume that $L$ is not normal in $G$. In particular, $L 
\not\supseteq G^{\prime}$, so $|L \cap G^{\prime}| \leq \frac{1}{2} 
|G^{\prime}|$. But since $G$ has class at most three and $Z \subseteq H$, we
see that $G^{\prime} \subseteq N$ and is abelian. Hence, by 
definition of $L$, $|G^{\prime} \cap T_{3,\alpha}| > \frac{1}{2}|G^{\prime}|$ and
so $G^{\prime} \subseteq T_{3,\alpha}$ since it is abelian. 
\par Now consider any $x \in T_{3,\alpha} \backslash H$ for which the coset $Hx$
is not exceptional. 
We shall show that $x \in N$,
which would imply that $N = G$, since there is at most one exceptional coset, 
contradicting our assumptions about $H$.
Let $h \in H$. Then, since $G^{\prime} \subseteq T_{3,\alpha}$, we have 
$[h,x]\alpha = [h,x]^3$. But also $[h,x]\alpha = [h\alpha,x\alpha] = 
[h^3,x^3]$ and, by Corollary \ref{cor:strange}, $[h^3,x^3] = [h,x]$. Thus
$[h,x]^2 = 1$ and another application of Corollary \ref{cor:strange}
implies that $h \in C_H(x^{-1}hx)$. But $C_H(x^{-1}hx) \supseteq C_H(x)$.
Thus, since $(H:C_H(x)) = 2$, we conclude that if $h \in H\backslash 
C_H(x)$ then $x^{-1}hx \in C_G(H)$, hence $x^{-1}hx \in H$ by Lemma 
\ref{lem:cosets}. Thus $x \in N$ as required, and this deals with 
{\sc Case 1}.
\\
\\
Now we turn to {\sc Case 2}. We have $|N \cap T_{3,\alpha}| \leq \frac{3}{4}
|N|$ by Lemma \ref{lem:strange}. On the other hand, Corollary 
\ref{cor:strange} and the fact that $G$ is nilpotent of class at most three
imply that every conjugate of $H$ lies in $N \cap T_{3,\alpha}$. To avoid
a contradiction we must have $(G:N) = 2$, thus $(G:H) = 4$. Write 
$G = H \sqcup Hx \sqcup Hy \sqcup Hz$, where $x,y,z \in T_{3,\alpha}$ and 
the cosets $Hx$ and $Hy$ are not exceptional. If $C_{H}(x) = C_{H}(y)$ then
$(H:Z) = 2$, a contradiction by Lemma \ref{lem:indextwo}. Otherwise, 
$(H:Z) = 4$ and, if $Hz$ is exceptional, then $Z = C_{H}(z)$ so that, in
particular, $z^2 \in Z$. Thus, by Corollary \ref{cor:strange}, the group 
$G/Z$, of order 16, has at least 8 involutions. In addition :
\par (i) $G/Z$ is non-abelian, since $G$ is not of class two,
\par (ii) $G/Z$ has a non-normal subgroup of order 4, namely $H/Z$,
\par (iii) $G/Z$ has no elements of order 8, since $G$ has no abelian subgroup 
of index two.
\\
These various restrictions serve to eliminate all possible structures for 
$G/Z$ (see \cite{TW}), a contradiction which completes the analysis of 
{\sc Case 2}.
\\
\\
Finally we turn to {\sc Case 3}. It is here that we at last will make use
of the induction hypothesis. If it were impossible to find 
$x_1, x_2 \in G\backslash N$ with $C_H(x_1) \neq C_H(x_2)$, then we'd have 
$(H:Z) = 2$, a contradiction by Lemma \ref{lem:indextwo}. So choose 
$x_1,x_2 \in G\backslash N$ with $C_H(x_1) \neq C_H(x_2)$ and such that
neither $Hx_1$ nor $Hx_2$ is exceptional, and pick any 
$h \in C_H(x_1) \backslash C_H(x_2)$. Consider the set 
\begin{eqnarray*}
S_h := \{g \in G : g^{-1}hg \in H \}.
\end{eqnarray*}
We have $N \subset S_h \subset G$, with all containments proper, since 
$x_1 \in S_h$ and $x_2 \not\in S_h$. But the fact that $G$ is nilpotent of 
class at most three, together with Lemma \ref{lem:center}, implies that 
$S_h$ is in fact a subgroup of $G$. Moreover it is $\alpha$-invariant, by
Corollary \ref{cor:strange}. Clearly $S_h$ satisfies the 
remaining hypotheses of Lemma \ref{lem:factorb} so, by minimality of
$G$, either $S_h$ has an abelian subgroup of index two or it is 
nilpotent class two. The former would imply that $N$ also contained an 
abelian subgroup of index 2, 
the latter that $H \lhd S_h$. Either way we have a contradiction, so 
the proof of Lemma \ref{lem:factorb} is complete.  
\end{proof}

It remains to prove Theorem \ref{thm:mainclass} under assumption 
\textbf{(C)}, that $G/Z$ contains an abelian subgroup of index 2. Let
this subgroup be $K/Z$ where $(G:K) = 2$. By Lemma \ref{lem:normal} we may 
assume a choice of $H$ which is not normal in $G$. Further we may 
assume that $(H:Z) > 2$, as otherwise, by Lemma \ref{lem:indextwo}, 
$G$ is clearly of type \textbf{II} in
Theorem \ref{thm:mainclass}. We consider two cases :
\\
\\
{\sc Case 1} : $K \supseteq H$.
\\
{\sc Case 2} : $(K:K\cap H) = 2$. 
\\
\\
If $K \supseteq H$ then $H \lhd K$ and so $K = N_G(H)$. Let 
\be
L := \{h \in H : g^{-1}hg \in H \; \forall \; g \in G \}.
\ee
Then $L$ is evidently a subgroup of $H$. We cannot have $L = H$ since 
otherwise $H \lhd G$. On the other hand, for any $x \in G\backslash 
K$ we have that $C_H(x) \subseteq L$. As there must be at least one 
non-exceptional coset of $H$ outside $K$, it follows that $(H:L) = 2$. 
If there is no exceptional coset, then clearly $L = Z$ and so 
$(H:Z) = 2$, a contradiction. Otherwise, notice that
$K$ is nilpotent class two and $\alpha$-invariant, being the normaliser
of $H$. Thus, by Lemma \ref{lem:normal}, it is of type 
\textbf{III} in Theorem \ref{thm:mainclass}. In particular, 
$(K:C_K(k)) \leq 2$ for every $k \in K$. Thus we'll still get the contradiction
that $L = Z$, unless $(K:H) = 2$ and $(H:Z) = 4$. So $|G/Z| = 16$ 
and $G$ is nilpotent. 
We can assume that 
\par (i) $G/Z$ is of class three, as otherwise $G$ would be of class at most
three and we could apply Lemma \ref{lem:factorb},
\par (ii) $G/Z$ has no elements of order 8, as otherwise $G$ would have an 
abelian subgroup of index 2, and thus clearly be of type \textbf{II}
in Theorem \ref{thm:mainclass}, since it is nilpotent,
\par (iii) $G/Z$ has a non-normal subgroup of order 4, namely $H/Z$. 
\\
These conditions eliminate all possible structures for $G/Z$ : see
\cite{TW}. We have dealt with {\sc Case 1}.
\\
\\
Finally, suppose $(K:K\cap H) = 2$. Let $H^{*} := K\cap H$. Clearly, 
$H^{*} \supseteq Z$ and $H^{*} \lhd G$. In fact $H^{*} = L$, the latter
defined as in (3.3). We can now argue as before, though note that there is 
an even easier approach : to avoid the contradiction that $(H:Z) = 2$ we'd need
to have $(G:N_G(H))=2$ and $H = N_G(H) \cap C_G(H^{*})$, which 
together yield the immediate contradiction that $H \lhd G$. 
\par This completes the proof of Theorem \ref{thm:mainclass}.     
  
\setcounter{equation}{0}
\section{Solvable Groups}\label{sec:theoremB}

In this section we further illustrate the effectiveness of the machinery
developed in Section 2 by proving

\begin{thm}\label{thm:solvable} Let $G$ be a finite group admitting an 
automorphism $\alpha$ for which $r_3(G,\alpha) > 4/15$. Then $G$ is 
solvable.
\end{thm}

The constant $4/15$ is best-possible, since $r_3(A_5,i) = 4/15$, where $i$
denotes the identity automorphism.
\\
\\
The proof of Theorem \ref{thm:solvable} is by induction on the group order.
Unsurprisingly, we shall have recourse to the 
classification of the finite simple groups in what follows, 
though the amount of 
information we draw on is quite limited and which we begin by summarising.

\begin{lem}\label{lem:abelian} Let $N$ be a non-abelian finite simple
group and $A$ an abelian subgroup of $N$ of maximum order. Then either
$|A|^3 < |N|$ or $N = L_2(q)$ for some prime power $q$, in which case
\begin{eqnarray*}
|N| = q(q^2-1)/2, \;\; |A| = q, \;\; {\hbox{if $q$ is odd}}, \\
|N| = q(q^2-1), \;\; |A| = q+1, \;\; {\hbox{if $q$ is even}}.
\end{eqnarray*}
In particular, $N$ has no abelian subgroup of index less than 12, and if
$N$ has an abelian subgroup of index less than 144, then 
$N = L_2(q)$ for some $q \in \{5,7,8,9,11,13\}$.
\end{lem}

\begin{proof}
The first assertion is the main result of \cite{V}. The second 
follows from a direct computation and the fact (see \cite{CCNPW}) that
the only non-abelian simple groups of order at most $\lfloor (143)^{3/2}
\rfloor = 1710$ are the groups $L_2(q)$ for $q \in \{5,7,8,9,11,13\}$. 
\end{proof}

\begin{lem}\label{lem:solvable} No non-abelian finite simple group possesses
a solvable subgroup of index less than 5. 
The only non-abelian finite simple groups
possessing a solvable subgroup of index at most 14 are the groups
$L_2(q)$, for $q \in \{5,7,8,9,11,13\}$, plus the group $L_3(3)$.
\end{lem}

\begin{proof} The first assertion follows from the solvability of $S_4$. 
For the second assertion, see \cite{CCNPW}.
\end{proof}

In the following table, $N$ is a non-abelian simple group, $A$ an abelian
subgroup of maximum order and $M$ a maximal subgroup of index at most 14. The
data and notation are taken from \cite{CCNPW}. 

\begin{tabular} {|c|c|c|c|c|c|} \hline
$N$ & $|N|$ & $(N:A)$ & $M$ & $(N:M)$ \\ \hline \hline
$L_2(5)$ & $2^2 \cdot 3 \cdot 5$ & $12$ & $A_4$ & $5$ \\ 
$\;$ & $\;$ & $\;$ & $D_{10}$ & $6$ \\ 
$\;$ & $\;$ & $\;$ & $S_3$ & $10$ \\ \hline
$L_2(7)$ & $2^3 \cdot 3 \cdot 7$ & $24$ & $S_4$ & $7$ \\
$\;$ & $\;$ & $\;$ & $\mathbb{Z}_7 \rtimes \mathbb{Z}_3$ & $8$ \\ \hline
$L_2(9)$ & $2^3 \cdot 3^2 \cdot 5$ & $40$ & $A_5$ & $6$ \\
$\;$ & $\;$ & $\;$ & $(\mathbb{Z}_3 \times \mathbb{Z}_3) \rtimes 
\mathbb{Z}_4$ & $10$ \\ \hline
$L_2(8)$ & $2^3 \cdot 3^2 \cdot 7$ & $56$ & 
$(\mathbb{Z}_2 \times \mathbb{Z}_2 \times \mathbb{Z}_2) \rtimes 
\mathbb{Z}_7$ & $9$ \\ \hline
$L_2(11)$ & $2^2 \cdot 3 \cdot 5 \cdot 11$ & $60$ & $A_5$ & $11$ \\ 
$\;$ & $\;$ & $\;$ & $\mathbb{Z}_{11} \rtimes \mathbb{Z}_5$ & $12$ \\ \hline
$L_2(13)$ & $2^2 \cdot 3 \cdot 7 \cdot 13$ & $84$ & $\mathbb{Z}_{13} 
\rtimes \mathbb{Z}_6$ & $14$ \\ \hline
$L_3(3)$ & $2^4 \cdot 3^3 \cdot 13$ & $432$ & 
$(\mathbb{Z}_3 \times \mathbb{Z}_3) \rtimes 
(\mathbb{Z}_2 \cdot S_4)$ & $13$ \\ \hline
\end{tabular}

$\;$ \\
$\;$ \\
From this table, we can also conclude the following :

\begin{lem}\label{lem:subgroups} Let $N$ be a non-abelian finite simple 
group. 
\par (i) If $N$ possesses a solvable subgroup $S$ such that $(N:S) \leq 14$
and $Z(S) \neq \{1\}$, then $N \cong L_2(5)$. 
\par (ii) If $N$ possesses a Sylow $2$-subgroup of index less than 45, then 
$N \cong L_2(5)$ or $L_2(7)$. 
\end{lem}

From now on, $G$ denotes a minimal
counterexample to Theorem \ref{thm:solvable} : 
our aim is to obtain a contradiction. We 
also fix a choice of $\alpha \in {\hbox{Aut}}(G)$ such that
$r_3(G,\alpha) > 4/15$. 

\begin{lem}\label{lem:simplesub} If $G$ has either an abelian subgroup of
index less than 144 or a solvable subgroup of index less than 25, then 
$G$ has a non-abelian characteristic simple subgroup $N$ with trivial
centraliser. Thus $G$ is isomorphic to a subgroup of 
Aut$(N)$. 
\end{lem}

\begin{proof}
By Lemmas \ref{lem:abelian} and \ref{lem:solvable}, any group 
satisfying either of the hypotheses of Lemma \ref{lem:simplesub} can possess
no subgroup of the form $N \times N$, where $N$ is a non-abelian simple 
group. In particular, either $G$ itself is simple, or it 
possesses a proper characteristic subgroup $N_1$. 
By Lemma \ref{lem:induced} and
the minimality of $G$, the factor group $G/N_1$ must be solvable. Thus 
$N_1$ must be insolvable. Repeating this argument, we see that either
$N_1$ is simple or possesses a proper characteristic subgroup $N_2$. Then 
$N_2$ is also characteristic in $G$, and so must be insolvable by 
Lemma \ref{lem:induced}. Iteration of the argument must terminate with a 
characteristic, non-abelian simple subgroup $N$ of $G$. Then $G$ possesses
a subgroup isomorphic to $N \times C_G(N)$. Our hypotheses on $G$ force
$C_G(N)$ to be solvable. But then $G/C_G(N)$ cannot be solvable, so 
$C_G(N) = \{1\}$ by minimality of $G$ and Lemma \ref{lem:induced}. This 
proves the lemma.
\end{proof}

Our idea to force a contradiction will be to use the information that 
$r_3(G,\alpha) > 4/15$ to produce a subgroup $S$ of $G$ which is either 
abelian of small index or solvable with non-trivial center and even smaller 
index. We then use the lemmas above to reduce the number of 
possibilities for $G$ to only a very few, which can be eliminated by 
direct computation. As in the previous section, we will work around a coset
decomposition with respect to a subgroup $H \subseteq T_{3,\alpha}$. However,
in this section $H$ will always be a cyclic group, rather than a 
subgroup of maximum order sitting inside $T_{3,\alpha}$. We will make more
forceful use of Proposition \ref{prop:numbertheory}, and to this end we
now introduce some more notation :
\\
\\
Let $n$ be a positive integer. We denote by $T(n)$ the maximum size of a 
subset of $\mathbb{Z}_n$, written additively, which contains 
no non-trivial solutions to either of the equations 
$a+b=2c$, $a+2b=3c$. We set $\tau_n := T(n)/n$. Roth's theorem
(see \cite{R}) implies that $\tau_n \rightarrow 0$ as $n \rightarrow 
\infty$. We have the following table of values :
\\
\\
\begin{tabular}{|c|c|c|c|} \hline
$n$ & $T(n)$ & $\tau_n$ \\ \hline \hline
$2$ & $1$ & $1/2$ \\ \hline
$4$ & $2$ & $1/2$ \\ \hline
$5$ & $2$ & $2/5$ \\ \hline
$7$ & $2$ & $2/7$ \\ \hline
$8$ & $2$ & $1/4$ \\ \hline
$10$ & $2$ & $1/5$ \\ \hline
$11$ & $2$ & $2/11$ \\ \hline
$13$ & $3$ & $3/13$ \\ \hline
$14$ & $3$ & $3/14$ \\ \hline
$16$ & $4$ & $1/4$ \\ \hline
$17$ & $4$ & $4/17$ \\ \hline
\end{tabular}

$\;$ \\
$\;$ \\
Also it is not difficult to verify that $\tau_n < 4/17$ for any $n > 17$. 
The proof of Theorem \ref{thm:solvable} will now be accomplished in a sequence
of steps, the goal of which is to progressively restrict the possible
orders of the elements in the subset $T_{3,\alpha}$ of our 
hypothetical counterexample $G$. At the end of this sequence of steps we
will be able to conclude that every element of $T_{3,\alpha}$ has order 
2 or 4. But then sending an element to its cube is the same as sending it to 
its inverse, so Theorem \ref{thm:solvable} follows from the analogous 
result in \cite{P}. 
\\
\\
\textbf{Step 1} : \emph{$T_{3,\alpha}$ contains no element of prime power order
$q$ where $q \geq 17$}. 
\\
\\
The arguments in this first step will provide a protoype for all remaining 
steps, so we present a careful reasoning here and later on become more
concise. Let $h \in T_{3,\alpha}$ be an element of prime-power order $q$ and 
$H := \; <h>$. We consider a decomposition of $G$ into right cosets of
$H$ and let $S$ be the subgroup of $G$ generated by all the cosets 
$Hx$ such that $Hx \cap T_{3,\alpha} \neq \phi$ and $C_H(x) \neq \{1\}$. 
Then $Z(S)$ is non-trivial, since $q$ is a prime power. Also
$S$ is $\alpha$-invariant. Let 
$(G:S) := r$ and $r_3(G,\alpha) := \xi$. By Proposition 
\ref{prop:numbertheory}, if $x \not\in S$ then $|Hx \cap T_{3,\alpha}|
\leq \tau_q|H|$. It follows that $r_3(G,\alpha) \leq \frac{\xi}{r} + 
\left(1-\frac{1}{r}\right)\tau_q$ and hence, since $r_3(G,\alpha) > 4/15$ we
must have 
\be
r < {\xi - \tau_q \over 4/15 - \tau_q}.
\ee
This is a non-trivial restriction whenever $\tau_q < 4/15$. Then in 
particular we must have $r_3(S,\alpha_S) > 4/15$ so, by minimality of 
$G$, either $S$ is solvable or $S = G$. But the latter 
contradicts minimality of $G$, by Lemma \ref{lem:induced}, since 
$Z(S) \neq \{1\}$. So we conclude that if $\tau_q < 4/15$ then 
$G$ contains a solvable subgroup $S$ with non-trivial center and of index 
bounded by (4.1). 
\par Now, as previously noted, if $q \geq 17$ then $\tau_q \leq 4/17$. 
Since $\xi \leq 1$ a priori, we then have, by (4.1), that 
$(G:S) \leq 24$. But, moreover, by Theorem \ref{thm:mainclass}, if 
$S$ is non-abelian then $\xi \leq 3/4$ and then (4.1) gives that
$(G:S) \leq 16$. So suppose $S$ is non-abelian. Since it is 
$\alpha$-invariant and solvable, so also is 
$S^{*} := {\hbox{Core}}_G(S)$ and $G/S^{*}$ is isomorphic to a subgroup of
$S_{16}$. But minimality of $G$ and Lemma \ref{lem:induced} force
$S^{*}$ to be trivial, hence $G$ itself is isomorphic to a subgroup of
$S_{16}$, contradicting the fact that $G$ possesses an element of 
prime power order $q \geq 17$. 
\par Thus $S$ must be abelian, i.e.: $G$ possesses an abelian subgroup of 
index at most 24. By Table 1 and Lemma \ref{lem:simplesub}, $G$ must then
be isomorphic to one of $L_2(5)$, $S_5$ and 
$L_2(7)$. One checks by direct calculation that none of these three 
groups possess an automorphism 
$\alpha$ such that $r_3(G,\alpha) > 4/15$. 
\\
\\
\textbf{Step 2} : $T_{3,\alpha}$ possesses no elements of order $13$. 
\\
\\
Suppose $h \in T_{3,\alpha}$ with $|h| = 13$. Let $H := \; <h>$ and consider
the subgroup $S$ defined in analogous manner to the previous step. 
Since $\tau_{13} = 3/13$, we obtain from (4.1) that either 
$S$ is abelian and $(G:S) \leq 21$ or $S$ is at least solvable with non-trivial
center and $(G:S) \leq 14$. The former possibility is dealt with as above. 
The latter implies, as above, that $G$ can be embedded in $S_{14}$. But then
the subgroup $H$ must be self-centralising in $G$, so $H = S$ and 
$|G| \leq 13 \times 14 = 182$. This leaves the same three possibilities for 
$G$, by Lemma \ref{lem:simplesub}, which have already been dealt with. 
\\
\\
\textbf{Step 3} : $T_{3,\alpha}$ possesses no elements of order $11$.
\\
\\
Suppose otherwise with $H = \; <h>$ and $|h| = 11$. Since 
$\tau_{11} = 2/11$, this time (4.1) yields $(G:S) \leq 9$ so that 
$G$ is embeddable in $S_9$, immediately contradicting the existence of any
element of order $11$ in $G$.
\\
\\
\textbf{Step 4} : $T_{3,\alpha}$ possesses no elements of order $16$.        
\\
\\
Suppose otherwise and let $H = \; <h>$ with $h \in T_{3,\alpha}$ and
$|h| = 16$. Since $\tau_{16} = 1/4 < 4/15$ we could proceed as before, but we
would be left with a greater number of possibilities for $G$ to check 
directly. So we modify our approach. One can check that, up to  
automorphisms, the only four-element subsets of $\mathbb{Z}_{16}$ that
avoid non-trivial solutions to both $a+b=2c$ and $a+2b=3c$ are 
\begin{eqnarray*}
\{0,1,4,5\}, \;\;\; \{0,1,4,13\}, \;\;\; \{0,1,5,12\}, \;\;\; \{0,1,12,13\}.
\end{eqnarray*}
The important point is that each of these sets contains either $4$ or $12$. 
Since the subset $\{4,12\}$ is characteristic in $\mathbb{Z}_{16}$, 
it now follows from Lemmas \ref{lem:abba} and \ref{lem:ap} that, for 
any $x \in T_{3,\alpha}$, 
\begin{eqnarray*}
|Hx \cap T_{3,\alpha}| + |Hx^{-1} \cap T_{3,\alpha}| \leq {7 \over 16} |H|.
\end{eqnarray*}
Note that this applies even if $Hx = Hx^{-1}$. As a result, when applying 
(4.1), we can replace $\tau_{16} = 1/4$ by the better constant 
$\tau^{\prime}_{16} = 7/32$, which will force the subgroup $S$ to be abelian
of index $16$ in $G$. Then the only two remaining possibilities for $G$
(namely $A_5$ and $S_5$) have already been considered. 
\\
\\
At this point we can state that any element of $T_{3,\alpha}$ must have order 
$2^i 5^k 7^k$, where $i \leq 3$ and $j,k \leq 1$. 
\\
\\
\textbf{Step 5} : $T_{3,\alpha}$ contains no elements of order $10, 14$ or $35$.
\\
\\
Let $p_1$, $p_2$ be distinct primes and let $h \in T_{3,\alpha}$ be an element
of order $p_1 p_2$. If $\tau_{p_{i}} \leq 1/2$ for $i = 1,2$, it is easy to
see that the argument in \textbf{Step 1} can be modified to produce a
subgroup $S$ of $G$ of index 
\begin{eqnarray*}
r < {\xi - \tau_{p_1 p_2} \over 4/15 - \tau_{p_1p_2}},
\end{eqnarray*}
which is still solvable, $\alpha$-invariant and, crucially, possessing 
non-trivial center. Indeed, $Z(S)$ will contain either $h^{p_1}$ or $h^{p_2}$. 
\par Now since, in fact, $\tau_n \leq 1/2$ for any $n > 1$, we can indeed apply
this argument. Of the numbers $\tau_{10}, \tau_{14}, \tau_{35}$, the largest
is $\tau_{14} = 3/14$ and this gives the worst bounds. We find then that
either $S$ is abelian of index at most $14$, or non-abelian of index 
at most $10$. The former option is subsumed by previous steps. For the 
latter, we deduce as before that $G$ is embeddable in $S_{10}$. But then
any cyclic subgroup of order $14$ is self-centralising, so 
$|G| \leq 14 \times 10 = 140$, and we're home and dry. 
\\
\\
At this point we can assume that any non-identity 
element of $T_{3,\alpha}$ has order
$2,4,5,7$ or $8$. 
\\
\\
\textbf{Step 6} : $T_{3,\alpha}$ contains no elements of order $5$. 
\\
\\
Since $\tau_{5} = 2/5 > 4/15$, we cannot apply the method of \textbf{Step 1}
directly. Let $H = \; <h>$ with $h \in T_{3,\alpha}$ and $|h| = 5$. I claim
that, if $T_{3,\alpha}$ contains no elements of order greater than $8$, then
for any $x \in T_{3,\alpha}$ with $[h,x] \neq 1$ we have in fact that 
$\mathscr{T}(H,x) = \{0\}$. Assuming this to be the case, we can then indeed
apply the method of \textbf{Step 1}, inserting $1/5$ in place of 
$2/5$ in eq. (4.1). This yields either an abelian subgroup $S$ of index
at most $12$, which is dealt with as before, or a solvable subgroup $S$
of index at most $8$ such that $H \subseteq Z(S)$. In the latter case, 
$G$ is embeddable in $S_8$ and so $|C_G(H)| \leq 30$. Thus 
$|G| \leq 240$ and we obtain no new possibilities for $G$ here either.
\par It thus remains to prove our claim. Suppose on the contrary that
$x \in T_{3,\alpha}$, $[h,x] \neq 1$ and $hx \in T_{3,\alpha}$. Then 
$(hx)\alpha^{2} = h^9 x^9 \in T_{3,\alpha}$. If $|x| \in \{2,4,8\}$ this implies
that $h^{-1}x \in T_{3,\alpha}$, a contradiction by Lemma \ref{lem:ap}. 
If $|x| = 5$ it implies that $h^{-1}x^{-1} \in T_{3,\alpha}$, 
contradicting Lemma \ref{lem:abba}. Finally, if $|x| = 7$ then instead
consider $(hx)\alpha^{4} = h^{81} x^{81}$. This is in $T_{3,\alpha}$, hence
so also is $hx^{4}$. But now each of $x,hx,hx^4$ and $hx^7$ is in 
$T_{3,\alpha}$, which contradicts Proposition \ref{prop:numbertheory} since
$(1,4,7)$ is an arithmetic progression. 

\newpage

So now every non-identity element of $T_{3,\alpha}$ may be assumed to have 
order $2,4,7$ or $8$. 
\\
\\
\textbf{Step 7} : $T_{3,\alpha}$ contains no elements of order $7$. 
\\
\\
We proceed as above. Assume $h \in T_{3,\alpha}$ with $|h| = 7$. Let 
$x \in T_{3,\alpha}$ and suppose that also $hx \in T_{3,\alpha}$. If
$|x| \in \{2,4,8\}$ then considering $(hx)\alpha^2 = h^9 x^9$ we have that
$h^2 x \in T_{3,\alpha}$, which forces $[h,x] = 1$ by 
Proposition \ref{prop:numbertheory}. If instead $|x| = 7$ then considering
$(hx)\alpha^3 = h^{27}x^{27}$ yields that $h^{-1}x^{-1} \in T_{3,\alpha}$, 
again forcing $[h,x] = 1$ by Lemma \ref{lem:abba}. 
\par Thus we can run through the method of \textbf{Step 1}, replacing 
$\tau_{7} = 2/7$ by the better constant $1/7$. We omit further details. 
\\
\\
We now come to the final step. Every element of $T_{3,\alpha}$ may be assumed 
to have order $2,4$ or $8$. In particular, every element of $T_{3,\alpha}$ has 
$2$-power order. 
\\
\\
\textbf{Step 8} : $T_{3,\alpha}$ contains no elements of order $8$. 
\\
\\
Suppose the contrary and consider the subgroup $S$ produced by the method of
\textbf{Step 1}. Since $\tau_8 = 1/4$, one easily 
checks that the following two possibilities arise :
\\
\par (i) $r_3(S,\alpha) \leq 1/2$ and $(G:S) \leq 14$. 
\par (ii) $r_3(S,\alpha) > 1/2$ and $(G:S) \leq 44$. 
\\
\\
First suppose (i) holds. Since $(G:S) < 25$ we can first apply Lemma 
\ref{lem:simplesub} to conclude that $G$ is isomorphic to a 
subgroup of Aut$(N)$, where $N$ is a non-abelian simple subgroup of $G$
isomorphic to one of the groups in Table 1. But now each of the groups in
this table has outer automorphism group of order at most 4 (see
\cite{CCNPW}). Since $Z(S)$ has non-trivial intersection with a 
cyclic group of order 8, this implies that the group $N \cap S$ also 
has a non-trivial center. But then Lemma \ref{lem:subgroups}(i) implies
that $N \cong A_5$, which we've already dealt with. 
\par Finally suppose (ii) holds. Since every element of $T_{3,\alpha}$ has 
$2$-power order, Theorem \ref{thm:mainclass} now implies that either $S$ is 
a $2$-group or possesses an abelian subgroup of index 2. 
In the former case, Lemmas \ref{lem:subgroups}(ii) and \ref{lem:simplesub} 
leave only
the same three possibilies for $G$ which were already encountered in 
\textbf{Step 1}, along with Aut$(L_2(7))$. This group is also eliminated
from consideration by direct calculation. In the latter case, if 
$S$ is not a 2-group then we must have $r_3(S,\alpha_S) \leq 2/3$, so
that (4.1) in fact yields $(G:S) \leq 24$. Thus $G$ possesses an 
abelian subgroup of index at most $48$, which leaves one more
possibility to rule out by direct calculation, namely $L_2(9) \cong A_6$. 
Having done so, the proof of Theorem \ref{thm:solvable} is 
complete. 

\setcounter{equation}{0}
\section{Conclusions}\label{sec:conclusions}

By establishing a connection to a well-studied problem in combinatorial 
number theory, Proposition \ref{prop:numbertheory} essentially forces
one of two alternatives on a finite group $G$ possessing an automorphism 
which cubes a large fraction of its elements : either a large fraction of 
pairs of elements in fact commute, or a large fraction of elements have small 
order. The connection to number theory may be interesting in its own right
for other reasons of which we are not aware. Also intruiging is whether 
Proposition \ref{prop:numbertheory} captures the full essence of this
connection, or whether there is more to be said. For example one could
ask to classify all minimal sets $\mathscr{S}$ of words in two letters
$a$ and $b$ such that, if an automorphism of a finite group $G$ cubes 
every element corresponding to a word in $\mathscr{S}$, then $a$ and $b$ must
commute. Are there any such sets other than those identified in Section 2 ? 
Also intruiging is whether Lemma \ref{lem:ap} holds in all infinite groups. 
We leave these matters for future investigation. 

\section*{Acknowledgement}

I thank Desmond MacHale for originally introducing me to these kinds of 
problems. Much of the material in Section 3 was completed a long time ago 
while a student of his, though I only spotted the connection to
number theory explicitly quite recently, and thus was able to prove 
Theorem \ref{thm:solvable} as well.

\ \\

\end{document}